\documentclass{article}
\usepackage{graphicx}
\usepackage{amsmath}
\usepackage{amsfonts}
\usepackage{amssymb}

\begin{document}

\begin{center}
{}\textbf{New proof of a Theorem on k-hypertournament losing scores}
\end{center}

$\bigskip$

\begin{center}
$^{1}$\textbf{S. Pirzada and \ }$^{2}$\textbf{Zhou Guofei}
\end{center}

\bigskip

$^{1}$Department of Mathematics, University of Kashmir, Srinagar-190006, India

Email: sdpirzada@yahoo.co.in \smallskip\bigskip

Department of Mathematics,Nanjing University, Nanjing, P.R.China

Email: gfzhou@nju.edu.cn\ 

\bigskip

\textbf{Abstract. }In this paper, we give a new and short proof of a Theorem
on k-hypertournament losing scores due to Zhou et al.[7].\bigskip

\begin{center}
\textbf{1. Introduction\bigskip}
\end{center}

An edge of a graph is a pair of vertices and an edge of a hypergraph is a
subset of the vertex set, consisting of atleast two vertices. An edge in a
hypergraph consisting of k vertices is called a k-edge, and a hypergraph all
of whose edges are k-edges is called a k-hypergraph.\bigskip

A k-hypertournament is a complete k-hypergraph with each k-edge endowed with
an orientation, that is, a linear arrangement of the vertices contained in the
hyperedge. In other words, given two non-negative integers n and k with $n\geq
k>1,$ a k-hypertournament on n vertices is a pair (V, A), where V is a set of
vertices with $\left|  V\right|  =n\ $and\ A is a set of k-tuples of vertices,
called arcs, such that any k-subset S of V, A contains exactly one of the k!
\ k-tuples whose entries belong to S.\ If n
$<$%
k, A = $\phi$\ and this type of \ hypertournament is called a
null-hypertournament. Clearly, a 2-hypertournament is simply a tournament.Let
e=(v$_{1},v_{2},...,v_{k})$ be an arc in a k-hypertournament H. Then
e(v$_{i},v_{j})$ represents the arc obtained from e by interchanging v$_{i}$
and v$_{j}.\bigskip$

The following result due to Landau [4] characterizes the score sequences in tournaments.\bigskip

\textbf{Theorem 1}. A sequence of non-negative integers [ s$_{1}$, s$_{2}%
$,..., s$_{n}$ ] in non-decreasing order is a score sequence of some
tournament if and only if

\begin{center}
$\sum_{i=1}^{j}s_{i}\geq\left(
\begin{array}
[c]{c}%
j\\
2
\end{array}
\right)  ,$\ \ \ 1 $\leq$ j $\leq$ n,
\end{center}

with equality when j = n.\bigskip

Now, there exist several proofs of Landau's theorem and a survey of these can
be found in Reid [5]. There are stronger inequalities on the scores in
tournaments which are due to Brualdi and Shen [1].\bigskip

Instead of scores of vertices in a tournament, Zhou et al.[7] considered
scores and losing scores of vertices in a k-hypertournament, and derived a
result analogous to Landau's theorem [4]. The score s(v$_{i}$) or s$_{i}$\ of
a vertex v$_{i}$\ is the number of arcs containing v$_{i}$\ and in which
v$_{i}$\ is not the last element, and the losing score r(v$_{i}$) or r$_{i}%
$\ of a vertex v$_{i}$\ is the number of arcs containing v$_{i}$\ and in which
v$_{i}$\ is the last element. The score sequence (losing score sequence) is
formed by listing the scores (losing scores) in non-decreasing order.

For two integers p and q, $\left(
\begin{array}
[c]{c}%
p\\
q
\end{array}
\right)  $ $=\frac{p!}{q!(p-q)!}$\ and $\left(
\begin{array}
[c]{c}%
p\\
q
\end{array}
\right)  $ $=0$\ \ \ if p
$<$%
q.\bigskip

The following characterization of losing score sequence in k-hypertournaments
are due to Zhou et al. [7].\bigskip

\textbf{Theorem 2.} Given two non-negative integers n and k with $n\geq
k>1,$\ a non-decreasing sequence R = [r$_{1}$, r$_{2}$ ,..., r$_{n}$] of
non-negative integers is a losing score sequence of some k-hypertournament if
and only if for each j,

\begin{center}
$\sum_{i=1}^{j}r_{i}\geq\left(
\begin{array}
[c]{c}%
j\\
k
\end{array}
\right)  ,$\ \ \ \ \ \ \ \ \ \ \ \ \ \ \ \ \ \ \ \ \ \ \ \ \ \ \ \ \ \ \ \ \ \ \ \ \ \ \ \ \ \ \ \ \ \ (1)
\end{center}

with equality when j = n.\bigskip

Koh and Ree [3 ] have given a different proof of Theorem 2. \bigskip Some
\ more results on scores of k-hypertournaments can be found in [2, 6]. The
following is the new and short proof of Theorem 2.\bigskip

\textbf{Proof. }The necessity part is obvious.

We prove sufficiency by contradiction. Assume all sequences of non-negative
integers in non-decreasing order of length fewer than n, satisfying conditions
(1) be the losing score sequences. Let n be the smallest length and r$_{1}$ be
the smallest possible with that choice of n such that R = [r$_{1}$, r$_{2}$
,..., r$_{n}$] is not a losing score sequence.

Consider two cases, (a) equality in (1) holds for some j
$<$%
n, and (b) each inequality in (1) is strict for all j
$<$%
n.

\textbf{Case (a).} Assume j (j
$<$%
\ n) is the smallest such that

\ \ \ \ \ \ \ \ \ \ \ \ \ \ \ \ \ \ \ \ \ \ \ \ \ \ \ \ \ \ \ \ \ \ \ \ \ \ $\sum
_{i=1}^{j}r_{i}\geq\left(
\begin{array}
[c]{c}%
j\\
k
\end{array}
\right)  .$

By the minimality of n, the sequence [r$_{1}$, r$_{2}$ ,..., r$_{j}$] is the
losing score sequence of some k-hypertournament H$_{1}$. Also,

\ \ \ \ \ \ \ \ \ \ $\sum_{i=1}^{m}\left[  r_{j+i}-\left(  \frac{1}{m}\right)
\sum_{i=1}^{k-1}\left(
\begin{array}
[c]{c}%
j\\
i
\end{array}
\right)  \left(
\begin{array}
[c]{c}%
n-j\\
k-i
\end{array}
\right)  \right]  $\ \bigskip

\ \ \ \ \ \ \ \ \ \ \ \ \ \ \ \ \ \ \ \ \ \ \ \ \ \ \ \ \ \ \ \ \ \ \ \ \ \ $=$%
\ $\sum_{i=1}^{m+j}r_{i}-\left(
\begin{array}
[c]{c}%
j\\
k
\end{array}
\right)  -\sum_{i=1}^{k-1}\left(
\begin{array}
[c]{c}%
j\\
i
\end{array}
\right)  \left(
\begin{array}
[c]{c}%
n-j\\
k-i
\end{array}
\right)  \bigskip$

\ \ \ \ \ \ \ \ \ \ \ \ \ \ \ \ \ \ \ \ \ \ \ \ \ \ \ \ \ \ \ \ \ \ \ \ \ \ $\geq
$\ $\left(
\begin{array}
[c]{c}%
m+j\\
k
\end{array}
\right)  -\left(
\begin{array}
[c]{c}%
j\\
k
\end{array}
\right)  -\sum_{i=1}^{k-1}\left(
\begin{array}
[c]{c}%
j\\
i
\end{array}
\right)  \left(
\begin{array}
[c]{c}%
n-j\\
k-i
\end{array}
\right)  \bigskip$

\ $\ \ \ \ \ \ \ \ \ \ \ \ \ \ \ \ \ \ \ \ \ \ \ \ \ \ \ \ \ \ \ \ \ \ \ \ \ =\left(
\begin{array}
[c]{c}%
m\\
k
\end{array}
\right)  .$

for each m, 1 $\leq$\ m $\leq$\ n-j, with equality when m = n-j.

Let $\frac{1}{m}\sum_{i=1}^{k-1}\left(
\begin{array}
[c]{c}%
j\\
i
\end{array}
\right)  \left(
\begin{array}
[c]{c}%
n-j\\
k-i
\end{array}
\right)  =\alpha.$ Therefore, by the minimality of n, the sequence
[r$_{k+1}-\alpha$, r$_{k+2}-\alpha$ ,..., r$_{n}-\alpha$] is the losing score
sequence of some k-hypertournament H$_{2}$.\ Taking disjoint union of H$_{1}$
and H$_{2}$, and adding all $m\alpha$\ arcs between H$_{1}$ and H$_{2}$ such
that each arc among $m\alpha$\ has the last entry in H$_{2}$ and each vertex
of H$_{2}$\ gets equal shares from these $m\alpha$\ last entries, we obtain a
k-hypertournament with losing score sequence R, which is a contradiction.

\textbf{Case (b).} Let each inequality in (1) is strict when j
$<$%
n, and in particular r$_{1}$
$>$%
0. Then the sequence [r$_{1}-1$, r$_{2}$ ,..., r$_{n}+1$] satisfies (1), and
therefore by minimality of r$_{1}$, is the losing score sequence of some
k-hypertournament H, a contradiction. Let x and y be the vertices respectively
with losing scores r$_{n}$+1 and r$_{1}$-1. If there is an arc e containing
both x and y with y as the last element in e, let e$^{/}$ = (x, y). Clearly,
(H-e)$\cup e^{/}$ is the k-hypertournament with losing score sequence R, again
a contradiction. If not, since r(x)
$>$%
r(y) there exist two arcs of the form e$_{1}$ = (w$_{1}$, w$_{2}$,...,
w$_{l-1}$, u, w$_{l}$,..., w$_{k-1}$) and e$_{2}$ = (w$_{1}^{/}$, w$_{2}^{/}%
$,..., w$_{k-1}^{/}$, v), where (w$_{1}^{/}$, w$_{2}^{/}$,..., w$_{k-1}^{/}$)
is a permutation of (w$_{1}$, w$_{2}$,..., w$_{k-1}$), x $\notin$\ \{w$_{1}$,
w$_{2}$,..., w$_{k-1}$\} and y $\notin$\ \{w$_{1}$, w$_{2}$,..., w$_{k-1}$\}.
Then, clearly R is the losing score sequence of the k-hypertournament
(H-(e$_{1}\cup$ e$_{2}$))$\cup$(e$_{1}^{/}\cup$ e$_{2}^{/}$), where e$_{1}%
^{/}=(u,w_{k-1})$, e$_{2}^{/}=(w_{t}^{/},v)$ and t is the integer with
$w_{t}^{/}=w_{k-1}.$ This again contradicts the hypothesis. Hence, the result follows.\bigskip

\ $\smallskip\ \ \ \ \ \ \ \ \ \ \ \ \ \ \ \ $

\begin{center}
\textbf{References\bigskip}\bigskip
\end{center}

[1] R.A.Brualdi and J.Shen, Landau's inequalities for tournament scores and a
short proof of a Theorem on transitive sub-tournaments, J. Graph Theory 38
(2001) 244-254. \ \ \ \ 

[2] \ Y. Koh and S. Ree, Score sequences of hypertournament matrices, J.Korea
Soc. Math. Educ. Ser. B: Pure and\ Appl. Math. 8 (2) (2001) 185-191.

[3] \ Y. Koh and S. Ree, On k-hypertournament matrices, Linear Algebra and its
Applications 373 (2003) 183-195.

[4] \ H. G. Landau, On dominance relations and the structure of animal
societies. III. The condition for a score structure, Bull. Math. Biophys.\ 15
( 1953 ) 143-148.

[5] K.B. Reid, Tournaments, scores, kings, generalizations and special topics,
Cong. Num.115 (1996) 171-211.

[6] \ C. Wang and G. Zhou, Note on the degree sequences of k-hypertournaments,
Discrete Mathematics, Preprint.

[7] \ G. Zhou, T. Yao and K. Zhang, On score sequences of k-hypertournaments,
European J. Combin. 21 (8)\ (2000) 993-1000.

\ \ \ \ \ \ \ \ \ \ \ \ \ \ \ \ \ 

\ \ \ \ 

\ \ \ \ \ \ \ \ \ \ \ \ \ \ \ \ \ \ \ \ \ \ \ \ \ \ \ \ \ \ \ \ \ \ 
\end{document}